\documentclass[10pt]{amsart}

\usepackage{amsfonts}
\usepackage{amssymb}
\usepackage{mathrsfs}
\usepackage{amsrefs}

\theoremstyle{plain}
\newtheorem{thm}{Theorem}[section]
\newtheorem{prop}[thm]{Proposition}
\newtheorem{cor}[thm]{Corollary}
\newtheorem{lemma}[thm]{Lemma}

\renewcommand{\mod}[1]{\left\lvert #1 \right\rvert}
\newcommand{\inp}[2]{\left\langle #1,#2 \right\rangle}
\newcommand{\norm}[1]{\left\| #1 \right\|}
\newcommand{\cl}[1]{\overline{#1}}
\newcommand{\bb}[1]{\mathbb{#1}}
\newcommand{\cal}[1]{\mathcal{#1}}

\DeclareMathOperator{\lat}{Lat}
\DeclareMathOperator{\Lat}{Lat}

\DeclareMathOperator{\linspan}{span}
\DeclareMathOperator{\lcm}{lcm}

\newcommand{\wk}{\text{weak$^\ast$}}
\DeclareMathOperator*{\mult}{mult}

\allowdisplaybreaks

\begin{document}

\title[Nevanlinna-Pick interpolation for $H^\infty_B$]{Nevanlinna-Pick interpolation for $\bb{C}+BH^\infty$}

\author[Mrinal Raghupathi]{Mrinal Raghupathi}
\address{Department of Mathematics, University of Houston \\ Houston, Texas 77204-3476, U.S.A.}
\email{mrinal@math.uh.edu}
\urladdr{http://www.math.uh.edu/~mrinal}
\thanks{This research was partially supported by the NSF grant DMS 0300128. This research was completed as part of my Ph.D. dissertation at the University of Houston.}

\begin{abstract}
We study the Nevanlinna-Pick problem for a class of subalgebras of $H^\infty$. This class includes algebras of analytic functions on embedded disks, the algebras of finite codimension in $H^\infty$ and the algebra of bounded analytic functions on a multiply connected domain. Our approach uses a distance formula that generalizes Sarason's~\cite{Sa} work. We also investigate the difference between scalar-valued and matrix-valued interpolation through the use of $C^\ast$-envelopes.
\end{abstract}

\maketitle

\section{Introduction}\label{intro}
Let us assume that we are given $n$ points $z_1,\ldots,z_n \in \bb{D}$, $n$ complex numbers $w_1,\ldots,w_n$ and a unital subalgebra $\cal{A}$ of $H^\infty$. We will say that a function $f\in \cal{A}$, \textit{interpolates} the values $z_1,\ldots,z_n$ to $w_1,\ldots,w_n$ if and only if $f(z_j)=w_j$. Such an $f$ will also be called an \textit{interpolating function}. We say that a function $f\in \cal{A}$ is a \textit{solution} to the interpolation problem if $f$ interpolates and $\norm{f}_\infty\leq 1$. 

For the algebra $H^\infty$, the Nevanlinna-Pick theorem gives an elegant criterion for the existence of a solution to the interpolation problem. The theorem states that a holomorphic function $f:\bb{D}\to\bb{D}$ interpolates $z_1,\ldots,z_n$ to $w_1,\ldots,w_n$ if and only if the Pick matrix 
\[\left[\frac{1-w_i\overline{w_j}}{1-z_i\overline{z_j}}\right]_{i,j=1}^n\]
is positive (semidefinite).

The operator theoretic approach to Nevanlinna-Pick interpolation has its roots in Sarason's work \cites{Sa}. Sarason used duality methods and Riesz factorization to prove a distance formula, which in turn implies the Nevanlinna-Pick theorem. Abrahamse~\cites{A} extended this approach to prove a Nevanlinna-Pick type interpolation theorem for the algebra $H^\infty(R)$ of bounded analytic functions on a multiply connected domain $R$. We will see in Section~\ref{examples} that the algebra $H^\infty(R)$ can be viewed as a subalgebra of $H^\infty$. These duality methods were used in \cites{DPRS} to prove an interpolation theorem for the codimension-1, $\wk$-closed, unital subalgebra of $H^\infty$ generated by $z^2$ and $z^3$. There is a clear connection in all three cases between the invariant subspaces for the algebra and the interpolation theorem.

It is too much to hope for a Nevanlinna-Pick theorem for all unital, $\wk$-closed subalgebras of $H^\infty$. In this paper we will study a certain class of subalgebras of $H^\infty$ that arise naturally. We will call these algebras the \textit{algebras with predual factorization}. We define these algebras in Section~\ref{distance}. 

In Section~\ref{examples} we will look at algebras of the form $H^\infty_B:=\mathbb{C}+BH^\infty$, where $B$ is an inner function and $\mathbb{C}$ denotes the span of the constant functions in $H^\infty$. This algebra is easily seen to be the unitization of $BH^\infty$, which is a $\wk$-closed ideal of $H^\infty$. This is the basic construction we will manipulate in order to obtain a large class of algebras with predual factorization. In section~\ref{invariant} we establish, for the algebra $H^\infty_B$, an analogue of the Helson-Lowdenslager theorem on invariant subspaces. We also establish the analogue of the Halmos-Lax theorem. We would like to point out that our invariant subspace theorem in Theorem~\ref{invsub} generalizes the result in~\cites{PS}. Our proof is also much more elementary and is in fact a simple consequence of the Helson-Lowdenslager theorem. In Section~\ref{distance} we compute the distance of an element $f\in L^\infty$ from the $\wk$-closed ideal $EH^\infty\cap \cal{A}$, where $E$ is the finite Blaschke product with zero set $z_1,\ldots,z_n$ and $\cal{A}$ is an algebra with predual factorization. As a consequence we obtain a Nevanlinna-Pick type theorem for $\cal{A}$. An important aspect of Nevanlinna-Pick interpolation is distinguishing between scalar-valued and matrix-valued interpolation theory. We will show in Section \ref{cstar} that the scalar-valued interpolation result obtained in Theorem~\ref{Hinftybinterp1} is not valid in the matrix-valued setting. Our work is inspired by, and considerably extends, the results in~\cites{DPRS}. 
\section{Notation and Examples}\label{examples}
Throughout this paper $L^p$, $1\leq p\leq \infty$, will denote the usual Lebesgue space of the circle. The subspace of elements whose negative Fourier coefficients are 0 will be denoted $H^p$ and we will freely use the identification of these spaces with spaces of analytic functions on the disk. We refer the reader to \cite{H} for the relevant background. We will regard $L^p$ as a normed space with its usual norm, but in dealing with $L^\infty$, and its subspaces, we will work with the $\wk$-topology that $L^\infty$ inherits as the dual of $L^1$. Given a non-empty set $\cal{S}\subseteq L^p$, we denote by $[\cal{S}]_p$ the smallest closed subspace of $L^p$ spanned by $\cal{S}$, when $p=2$ we will denote this as $[\cal{S}]$ and when $p=\infty$ we use  $[\cal{S}]_\infty$ to denote the $\wk$-closed subspace spanned by $\cal{S}$.

Let $g\in H^\infty$ and let $H^p_g$ denote the space $[\bb{C}\cdot 1+g H^p]_p$. If $u$ is an outer function, then $[uH^p]_p=H^p$. If we factor $g$ into its inner factor $B$ and outer factor $u$, then $H^p_g=[\bb{C}\cdot 1+g H^p]_p=[\bb{C}\cdot 1+BH^p]_p=H^p_B$. Therefore it is enough to consider inner functions. 

Our primary interest will be in the spaces $H^p_B$ for a Blaschke product $B$ and for $p=1,2,\infty$. We denote by $\phi_a$ the elementary M\"{o}bius transformation of the disk given by $\phi_a(z)=\dfrac{a-z}{1-\bar{a}z}$, for $a\not =0$, and $\phi(z)=z$, when $a=0$. We will write the Blaschke product $B$ as
\[B=\prod_{j\in J} \dfrac{\mod{\alpha_j}}{\alpha_j}\phi_{\alpha_j}^{m_j},\] 
where $J$ is either a finite or countably infinite set, $\{\alpha_j\}$ are distinct and $m_j\geq 1$. The normalizing factor $\dfrac{\mod{\alpha_j}}{\alpha_j}$ is introduced to ensure convergence and is defined to be $1$, if $\alpha_j=0$. We will assume throughout that $B$ has at least 2 zeros, i.e., $\sum_{j=in J} m_j\geq 2$. We point out that a function $f\in H^\infty_B$ if and only if it satisfies the following twoc onstraints:
\begin{enumerate}\label{constraints}
\item $f(\alpha_i)=f(\alpha_j)$, $i,j\in J$.
\item If $m_j\geq 2$, then $f^{(i)}(\alpha_j)=0$ for $i=1,\ldots,m_j-1$.
\end{enumerate}
We will use $K$ to denote the Szeg\"o kernel and $k_z$ to denote the Szeg\"o kernel at the point $z$, i.e., the element of $H^2$ such that $f(z)=\inp{f}{k_z}$ for all $f\in H^2$. Note that the Szeg\"o kernel is actually a bounded analytic function on $\bb{D}$ and so is in $H^\infty$. We will often abuse notation by using the letter $z$ to represent a complex variable, the identity map on $\bb{D}$ and the identity map on $\bb{T}$.

We now provide some examples from the literature that motivate these definitions. 

First, consider the subalgebra $H^\infty_1$ of functions for which $f'(0)=0$, this algebra has been studied in \cites{DPRS}. Clearly, $H^\infty_1=\bb{C}+z^2H^\infty$ and corresponds to the case where $B$ is the Blaschke product $z^2$. The algebra $H^\infty_1$ is generated by $z^2$ and $z^3$.

Second, consider the algebra of functions that are equal at the points $a,b\in \bb{D}$ where $a$ is different from $b$. Note that if $f\in H^\infty$ and $f(a)=f(b)$, then $f-f(a)$ vanishes at both $a$ and $b$. Hence, $f-f(a)=Bh$, where $h\in H^\infty$ and $B$ is the Blaschke product for the points $a,b$. Hence, this algebra is $\bb{C}+BH^\infty$. This example can be extended to any Blaschke sequence. Interpolation questions for algebras of this type where the zeros of the Blaschke product are finite in number, and of multiplicity one, were studied in~\cites{So}. In~\cites{So} the focus was the matrix-valued theory and an additional condition was imposed that reduced the result to the classical situation. When we study the matrix-valued interpolation problem in Section~\ref{cstar} we will need to impose a similar restriction.

The two examples just mentioned are special cases of the cusp algebras studied by Agler and McCarthy~\cites{AM1,AM2}. Cusp algebras are of finite codimension in $H^\infty$. We do not require $H^\infty_B$ to be of finite codimension. In dealing with cusp algebras and embedded disks we must consider not only the algebra $H^\infty_B$ but also consider finite intersections of such algebras. The interpolation theorem that we prove in Section~\ref{distance} applies to both the infinite codimension case as well as an infinite intersection of algebras of the form $H^\infty_B$.

Our third example provides some additional motivation as to why one may wish to study $H^\infty_B$. If $G$ is a group of conformal automorphisms of the disk $\bb{D}$, then we have a natural action of the group $G$ on $H^\infty$ given by $f\mapsto f\circ \alpha^{-1}$, where $\alpha\in G$. The fixed point algebra is denoted $H^\infty_G$. In the case where the group is a Fuchsian group there is a natural Riemann surface structure on $R=\bb{D}/G$ and $H^\infty_G$ is isometrically isomorphic to the algebra of bounded analytic functions on $R$. If we assume for a moment that $H^\infty_G$ is non-trivial, then it must be the case that $\sum_{\alpha\in G}(1-\mod{\alpha(\zeta)})<\infty$ for any $\zeta\in \bb{D}$. Let $B_\zeta$ be the Blaschke product whose zero set is the points of the orbit of $\zeta$ under $G$. The Blaschke product $B_\zeta$ converges for each $\zeta \in \bb{D}$. From the constraints on page~\pageref{constraints} we see that $\bb{C}+B_\zeta H^\infty$ is the set of functions in $H^\infty$ that are constant on the orbit of $\zeta$, and so $H^\infty_G=\bigcap_{\zeta \in \bb{D}} H^\infty_{B_{\zeta}}$. 

While the interpolation results presented in Section~\ref{distance} do apply to the algebra $H^\infty_G$, the result is not as refined as Abrahamse's theorem~\cite{A}. A result that generalizes Abrahamse's theorem to the algebra $H^\infty_G$ can be found in~\cite{MRthesis}. 

\section{Invariant Subspaces}\label{invariant}

In this section we examine the invariant subspaces for the algebra $H^\infty_B$. 

Let $\cal{H}$ be a Hilbert space, let $\cal{A}$ be a subalgebra of $B(\cal{H})$. We say that a subspace $\cal{M}\subseteq \cal{H}$ is invariant for $\cal{A}$ if and only if $A(\cal{M})\subseteq \cal{M}$ for all $A\in\cal{A}$. 

We first make some general comments about models for invariant subspaces for subalgebras of $H^\infty$. The algebra $H^\infty$ has an isometric representation on $B(H^2)$. Suppose that $\cal{A}$ is a subalgebra of $H^\infty$. For the purposes of invariant subspaces we may assume that $\cal{A}$ is unital and $\wk$-closed, since the lattice of invariant subspaces is unaffected by unitization and $\wk$-closure. If $\cal{M}\subseteq H^2$ is invariant for $\cal{A}$, then Beurling's theorem tells us that $\cal{M}=[\cal{A}\cal{M}]\subseteq [H^\infty \cal{M}]=\phi H^2$ for some inner function $\phi$. Hence, $\cal{M}=\phi \cal{N}$ for some subspace $N$. It is the subspace $\cal{N}$ and not the inner factor $\phi$ which is relevant. In some sense the class of subspaces $\cal{N}$ that can arise are the models for invariant subspaces. 

A natural class of invariant subspaces for an algebra $\cal{A}$ are the cyclic subspaces. If $f\in H^2$, then we know that it has an inner-outer factorization $f=\phi u$ which is unique up to multiplication by a unimodular scalar. The cyclic subspace $\cal{M}=[\cal{A}f]=\phi[\cal{A}u]$. We see that the cyclic subspaces generated by $\cal{A}$ and an outer function $u$ are the models for all cyclic subspaces for $\cal{A}$. This will be relevant to our discussion about interpolation in Section~\ref{interp}.

Given a subalgebra $\cal{A}\subseteq H^\infty$ we say that a subspace $\cal{M}\subseteq L^p$ is invariant for $\cal{A}$ if and only if $f g\in \cal{M}$, whenever $f\in \cal{A}$ and $g\in\cal{M}$. We will assume, unless stated otherwise, that the term subspace means closed, non-trivial subspace.

We begin our study with a look at the subspaces of $H^2$ that are invariant for $H^\infty_B$. Beurling's theorem proves that the shift invariant subspaces of $H^2$,  i.e., the subspaces of $H^2$ invariant for the algebra $H^\infty$,  are of the form $\phi H^2$, where $\phi$ is inner. The Helson-Lowdenslager theorem is an extension of Beurling's theorem. Helson and Lowdenslager proved that the simply shift invariant subspaces of $L^2$ are of the form $\phi H^2$ where $\phi$ is unimodular, i.e., $\mod{\phi} = 1$ a.e.~on $\bb{T}$. Our first result is the analogue of the Helson-Lowdenslager theorem. 

\begin{thm}\label{invsub}
Let $B$ be an inner function and let $\cal{M}$ be a subspace of $L^p$ which is invariant under $H^\infty_B$. Either there exists a measurable set $E$ such that $\cal{M}=\chi_EL^p$ or there exists a unimodular function $\phi$ such that $\phi B H^p \subseteq \cal{M} \subseteq \phi H^p$. In particular, if $p=2$, then there exists a subspace $W\subseteq H^2\ominus BH^2$ such that $\cal{M}=\phi(W\oplus BH^2)$.
\end{thm}

\begin{proof}
The space $[B H^\infty \cal{M}]_p$ is a shift invariant subspace of $L^p$ and since $B$ is an inner function $[B H^\infty \cal{M}]_p=B[H^\infty \cal{M}]_p$. By the invariant subspace theorem for  $H^\infty$, either $[H^\infty \cal{M}]_p=\chi_E L^p$ for some measurable subset $E$ of the circle or $[H^\infty \cal{M}]_p=\phi H^p$ for some unimodular function $\phi$. In the former case 
\[\cal{M}\supseteq B[H^\infty \cal{M}]_p=B\chi_E L^p=\chi_E L^p\supseteq \cal{M}.\]
In the latter case we see that 
\[\phi BH^p = B[H^\infty \cal{M}]_p \subseteq \cal{M}\subseteq [H^\infty \cal{M}]=\phi H^p.\] 
When $p=2$, since $B H^2 \subseteq \overline{\phi}\cal{M}\subseteq H^2$, we see that $\cal{M}=\phi(W\oplus B H^2)$ where $W\subseteq H^2\ominus B H^2$. 
\end{proof}

As a corollary we obtain:

\begin{cor}[\cite{DPRS}*{Theorem 2.1}]
Let $H^\infty_1$ denote the algebra of functions in $H^\infty$ such that $f'(0)=0$. A subspace $\cal{M}$ of $L^2$ is invariant under $H^\infty_1$, but not invariant under $H^\infty$, if and only if there exists an inner function $\phi$, scalars $\alpha,\beta\in \bb{C}$ with $\mod{\alpha}^2+\mod{\beta}^2=1$, with $\alpha\not=0$, such that $\cal{M}=\phi([\alpha+\beta z]\oplus z^2H^2) $.
\end{cor}

\begin{proof}
From the previous result we see that $\cal{M}=\phi(W\oplus z^2H^2)$ where $W\subseteq H^2\ominus z^2H^2=\linspan\{1,z\}$. Since $\cal{M}$ is not invariant under $H^\infty$ see that $W$ is one-dimensional and that $\alpha\not=0$.
\end{proof}

We will identify inner functions that differ only by a constant factor of modulus 1.  If $\cal{S}\not=\{0\}$ is a subset of $H^2$, then Beurling's theorem tells us that $[H^\infty \cal{S}]=\phi H^2$ for some inner function $\phi_\cal{S}$. The inner function $\phi_\cal{S}$ is called the \index{Inner divisor}\textit{inner divisor} of $\cal{S}$. If $\cal{S}_1$ and $\cal{S}_2$ are two subsets of $H^2$, then we define their \index{Greatest common divisor}\textit{greatest common divisor} $\gcd(\cal{S}_1,\cal{S}_2)$ to be the inner divisor of $[H^\infty(\cal{S}_1\cup \cal{S}_2)]$ and the \index{Least common multiple}\textit{least common multiple} $\lcm(\cal{S}_1,\cal{S}_2)$ to be the inner divisor of $[H^\infty \cal{S}_1]\cap [H^\infty \cal{S}_2]$. For a function $f$ the inner divisor of $\{f\}$ is clearly the inner factor of $f$. For functions $f_1,f_2\in H^2$ we define $\gcd(f_1,f_2):=\gcd(\{f_1\},\{f_2\})$ and $\lcm(f_1,f_2):=\lcm(\{f_1\},\{f_2\})$. For a more detailed description of these operations we refer the reader to~\cite{bercovici}.

Let $\cal{A}\subseteq B(\cal{H})$ be an operator algebra. Associated to this operator algebra is its \index{Lattice of invariant subspaces}\textit{lattice of invariant subspaces}, which is defined as the set of subspaces of $\cal{H}$ that are invariant for $\cal{A}$. We will denote the lattice of non-trivial, invariant subspaces of $\cal{A}$ by $\lat(\cal{A})$. 

An important consequence of  Beurling's theorem is that it allows a complete description of the lattice of invariant subspaces for $H^\infty$. Two shift invariant subspaces $\phi H^2$ and $\psi H^2$ are the equal if and only if $\phi = \lambda \psi$ for a unimodular constant $\lambda$. Since we have chosen to identify inner functions that differ only by a constant, we see that that the shift invariant subspaces of $H^2$ are parametrized by inner functions. There is a natural ordering of inner functions. If $\phi,\psi$ are inner functions, then we say that $\phi\leq \psi$ if and only if there exists an inner function $\theta$ such that $\phi \theta  = \psi$. This makes the set of inner functions a lattice with meet and join given by 
\[
\phi \wedge \psi = \gcd(\phi,\psi), \phi \vee \psi = \lcm(\phi,\psi).
\]
In this ordering the inner function $1$ is the least element of the lattice and the lattice has no upper bound. The map $\phi\mapsto \phi H^2$ is a bijection between the lattice of inner functions and the the lattice of non-trivial, invariant subspaces for $H^\infty$. This identification is a lattice anti-isomorphism, i.e., order reversing isomorphism, taking meets to joins and joins to meets.

For the lattice $\Lat(H^\infty_B)$ the situation is different. There are two parameters that determine an invariant subspace $\cal{M}\in \Lat(H^\infty_B)$, an inner function $\phi$ and a subspace $W\subseteq H^2\ominus BH^2$. However, the subspace $\cal{M}$ does not uniquely determine $\phi$ and $W$. Conversely, different choices of $\phi$ and $W$ can sometimes give rise to the same subspace. A simple example is obtained by setting $B=z^2$, in which case 
\begin{equation}\label{subspace1}
zH^2=z\left([1,z]\oplus z^2H^2\right) = [z]\oplus z^2H^2.
\end{equation} 
In general, if $\cal{M}=\phi(W\oplus BH^2)$, then the subspace $W=\cl{\phi}\cal{M}\ominus BH^2$. It is always possible to make a canonical choice of inner function and subspace $W$. The canonical choice is to set the inner function equal to $\phi_\cal{M}$, the inner divisor of $\cal{M}$, and to let $W_\cal{M}=\cl{\phi_\cal{M}}\cal{M}\ominus BH^2$.

We now describe the extent to which the decomposition of the subspace $\cal{M}$ into the form $\phi(W\oplus BH^2)$ fails to be unique. It is useful to keep in mind the rather trivial example in~\eqref{subspace1}. Note that in addition to being invariant for $H^\infty_{z^2}$, the subspace $zH^2$ is also shift invariant.

\begin{prop}\label{lathinftyb}
Let $\cal{M}\in \Lat(H^\infty_B)$, let $\phi_\cal{M}$ be the inner divisor of $\cal{M}$ and let $W_\cal{M}=\cl{\phi_\cal{M}}\cal{M}\ominus BH^2$. Let $\psi$ be inner and $V$ be a subspace of $H^2\ominus BH^2$ such that $\cal{M}=\psi(V\oplus BH^2)$. The following are true:
\begin{enumerate}
\item The inner function $\gcd(W_\cal{M},B)=1$.
\item  The inner function $\psi\leq\phi_\cal{M}$ and 
\begin{equation}
\label{subspace2}\phi_\cal{M}W_\cal{M}=\psi V\oplus  B(\psi H^2\ominus \phi_\cal{M} H^2).
\end{equation}
\item If $\theta$ is such that $\psi\theta = \phi_\cal{M}$, then $\theta = gcd(B,V)$.
\item We have $\phi_\cal{M}=\psi$ if and only if $W_\cal{M}=V$.
\item If $\cal{M}\not\in \Lat(H^\infty_C)$ for all $C<B$ , then $\psi = \phi_\cal{M}$ and $V=W_\cal{M}$. \end{enumerate}
\end{prop}

\begin{proof} \rule{1in}{0mm}
\begin{enumerate}
\item Note that $\phi_\cal{M}\gcd(W_\cal{M},B)$ is an inner function that divides $\cal{M}$. Since $\phi_\cal{M}$ is the inner divisor of $\cal{M}$ we get $\gcd(W_\cal{M},B)=1$.
\item Since $\phi$ is the inner divisor of $\cal{M}$, it follows that $\psi|\phi$. Let $\theta$ be the inner function such that $\psi \theta = \phi$. We have,
\begin{equation*}
\psi\theta (W_\cal{M}\oplus BH^2) = \phi_\cal{M}(W_\cal{M}\oplus BH^2) = \psi (V\oplus BH^2).
\end{equation*}
It follows that
\begin{equation*}
\theta W_\cal{M}\oplus \theta BH^2  = V\oplus BH^2  = V\oplus B(H^2\ominus \theta H^2)\oplus B\theta H^2.
\end{equation*}
Hence, 
\begin{equation}\label{subspace3}
\theta W_\cal{M}=V\oplus B(H^2\ominus \theta H^2).
\end{equation}
Multiplying by $\psi$ gives~\eqref{subspace2}. 

\item From~\eqref{subspace3} we see that $\theta$ divides both $V$ and $B$ and so $\theta\leq \gcd(B,V)$. From~\eqref{subspace3} we get that $\gcd(B,V)|\theta W_\cal{M}$. Since $\gcd(W_\cal{M},B)=1$ it must be the case that $\gcd(B,V)\leq \theta$. Hence, $\theta = \gcd(B,V)$.

\item The conditions $\psi=\phi_\cal{M}$ and $W_\cal{M}=V$ are equivalent. If $\psi=\phi_\cal{M}$, then~\eqref{subspace3} shows that $W_\cal{M}=V$. Conversely, if $W_\cal{M}=V$, then~\eqref{subspace3} shows that $\theta W_\cal{M}\supseteq W_\cal{M}$. If $w\in W_\cal{M}\subseteq \theta W_\cal{M}$, then there exists $w_1\in W_\cal{M}$ such that $w=\theta w_1$. Repeating the argument we find that there exists $w_n\in W_\cal{M}$ such that $\theta^nw_n = w$. If $\theta\not =1$, then the equation $\theta^nw_n = w$ for all $n\geq 0$, contradicts the fact that $\theta$ cannot divide $w$ with infinite multiplicity. Hence, $\theta=1$ and $\phi_\cal{M} =\psi\theta = \psi$.

\item If $\theta \not =1$, then 
\begin{equation*}
\cal{M} = \psi(V\oplus BH^2) = \psi\theta(X \oplus CH^2),
\end{equation*}
where $C<B$. Hence, $\cal{M}\in\Lat(H^\infty_C)$. 
\end{enumerate}
\end{proof}

Proposition~\ref{lathinftyb} indicates that the lattice of invariant subspaces for $H^\infty_B$ is more complicated than the lattice of shift invariant subspaces. Although the canonical choice of inner divisor seems natural, this choice does not behave as expected with respect to the lattice operations. We do not, as yet, have at our disposal a useful way to describe the lattice operation in $\Lat(H^\infty_B)$. For illustrative purposes we examine what happens to the inner divisor when we take meets and joins of elements in $\Lat(H^\infty_B)$. Note that $\Lat(H^\infty)$ is a sublattice of $\Lat(H^\infty_B)$. Any good description of $\Lat(H^\infty_B)$ would have to take into account the fact that $\Lat(H^\infty)$ is the lattice of inner functions.

Let $\cal{M}=\phi_\cal{M}(W_\cal{M}\oplus BH^2),\cal{N}=\phi_\cal{N}(W_\cal{N}\oplus BH^2)\in \Lat(H^\infty_B)$, where $\phi_\cal{M}$ and $\phi_\cal{N}$ are the inner divisors of $\cal{M}$ and $\cal{N}$ respectively. Let $\cal{X}=\cal{M}\cap\cal{N}$. We have,
\begin{align*}
B\lcm(\phi_\cal{M},\phi_\cal{N})H^2&\subseteq \lcm(B\phi_\cal{M},B\phi_\cal{N})H^2=(B\phi_\cal{M} H^2)\cap (B\phi_\cal{N} H^2)\\
&\subseteq \cal{M}\cap \cal{N} =\cal{X}\\
&\subseteq (\phi_\cal{M} H^2)\cap(\phi_\cal{N} H^2)=\lcm(\phi_\cal{M},\phi_\cal{N})H^2.
\end{align*}
Hence,  $\phi_\cal{X}$ satisfies 
\begin{equation*}
\lcm(\phi_\cal{M},\phi_\cal{N})\leq \phi_\cal{X} \leq B\lcm(\phi_\cal{M},\phi_\cal{N}).
\end{equation*}

These are the best general bounds we have. If we consider the case where $\phi_\cal{M}=\phi_\cal{N}=1$, then we see that 
\begin{align*}
\cal{X}&=\cal{M}\cap \cal{N}\\
&= (W_\cal{M}\cap W_\cal{N})\oplus BH^2\\
&= \gcd(W_\cal{M}\cap W_\cal{N},B)(W_\cal{X}\oplus BH^2) 
\end{align*}
If $W_\cal{M}\cap W_\cal{N}=\{0\}$, then the inner divisor $\phi_\cal{X}=B$. However, if $W_1\cap W_2$ is non-trivial the situation can be different. Let $B=z^5$, let $\cal{M}=[1+z^2,z^3]\oplus z^5H^2$ and let $\cal{N}=[1-z^2,z^3]\oplus z^5H^2$. It is straightforward to check that the inner divisor $\phi_\cal{X}$ of the intersection $\cal{X}=\cal{M}\cap\cal{N}$ is divisible by $z^3$. Since the functions $1+z^2$ and $1-z^2$ are outer we see that $\phi_\cal{M}=\phi_\cal{N}=1$. Note that $\gcd(\phi_\cal{M},\phi_\cal{N})=1<\phi_\cal{X}<z^5=B$.

If we consider the join of two subspaces $\cal{Y}=\cal{M}\vee \cal{N}$, then we have $\phi_\cal{Y}=\gcd(\phi_\cal{M},\phi_\cal{N})$. The inequality $\gcd(\phi_\cal{M},\phi_\cal{N})\leq \phi_\cal{Y}$ follows from 
\begin{equation*}
\cal{M}\vee\cal{N}\subseteq (\phi_\cal{M} H^2) \vee (\phi_\cal{N} H^2) = \gcd(\phi_\cal{M},\phi_\cal{N})H^2.
\end{equation*}
Since $\phi_\cal{Y} |\cal{Y}$, we have $\phi_\cal{Y} |\cal{M}$ and $\phi_\cal{Y} |\cal{N}$. Hence, $\phi_\cal{Y} |\gcd(W,B)\phi_\cal{M}=\phi_\cal{M}$ and $\phi_\cal{Y}|\gcd(V,B)\phi_\cal{N}=\phi_\cal{N}$. Therefore, $\phi_\cal{Y}\leq \gcd(\phi_\cal{M},\phi_\cal{N})$, which implies $\phi_\cal{Y}=\gcd(\phi_\cal{M},\phi_\cal{N})$.

We now turn our attention to the vector-valued invariant subspaces for $H^\infty_B$. It is not difficult to extend Theorem \ref{invsub} to the vector-valued setting. If $\cal{H}$ is a separable Hilbert space, then we denote by $H^2_\cal{H}$ the $\cal{H}$-valued Hardy space. The natural action of $H^\infty$ on $H^2_\cal{H}$ is given by $(fh)(z)=f(z)h(z)$ and this makes $H^2_\cal{H}$ a module over $H^\infty$. This action obviously restricts to $H^\infty_B$ and we say that a subspace $\cal{M}$ of $H^2_\cal{H}$ is invariant for $H^\infty_B$ if and only if $H^\infty_B \cal{M}\subseteq \cal{M}$. We denote by $H^\infty_{B(\cal{H})}$ the set of $B(\cal{H})$-valued bounded analytic functions. An element of $H^\infty_{B(\cal{H})}$ is called rigid if $\Phi(e^{i\theta})$ is a partial isometry a.e.~on $\bb{T}$. A subspace $\cal{M}$ is invariant under $H^\infty$ if and only if there exists a rigid function $\Phi\in H^\infty_{B(\cal{H})}$ such that $\cal{M}=\Phi H^2_\cal{H}$. The proof of the scalar case carries through with the obvious modifications to give the following result. 

\begin{thm}
If $\cal{M}$ be a closed subspace of $H^2_\cal{H}$ which is invariant for $H^\infty_B$, then there exists a rigid function $\Phi \in H^\infty_{B(H)}$ and a subspace $V\subseteq H^2_\cal{H}\ominus BH^2_\cal{H}$ such that $\cal{M}=\Phi(V\oplus BH^2_\cal{H})$.
\end{thm} 

\begin{proof}
Let $\cal{M}\subseteq H^2_\cal{H}$ be an invariant subspace for $H^\infty_B$. As in the proof of Theorem~\ref{invsub} we form the shift invariant subspace $[H^\infty \cal{M}]\subseteq H^2_\cal{H}$. By the invariant subspace theorem in~\cite{helsoninvariant} we can write $[H^\infty \cal{M}]=\Phi H^2_\cal{H}$ for a rigid function $\Phi$. Now, $\cal{M} \supseteq [BH^\infty \cal{M}] = B[H^\infty \cal{M}]=B\Phi H^2_\cal{H}$ and so $B\Phi H^2_\cal{H}\subseteq \cal{M}\subseteq \Phi H^2_\cal{H}$. It follows that $\cal{M} = \cal{W}\oplus B\Phi H^2_\cal{H}$, where $\cal{W}\subseteq \Phi H^2_\cal{H}\ominus B\Phi H^2_\cal{H}$. If $w\in \cal{W}$, then $w=\Phi f$ for some $f\in H^2_\cal{H}\ominus BH^2_\cal{H}$. Choosing $V$ to be the subspace of elements $f\in H^2_\cal{H}\ominus BH^2_\cal{H}$ such that $\Phi f\in \cal{W}$ completes the proof.
\end{proof}

Keeping in mind our comments about models for invariant subspaces we see that $H^2$ serves as a model for subspaces of $H^2$ invariant for $H^\infty$. For the algebra $H^\infty_B$ the situation is more complicated. Even in the simplest cases, for example $H^\infty_{z^2}$, there can be infinitely many models. 

In the vector-valued case, $H^2_\cal{H}$ which is a direct sum of $H^2$ spaces, forms the only model for an invariant subspace of $H^\infty$. The models for invariant subspaces of $H^\infty_B$ are parametrized by subspaces $V\subseteq H^2_\cal{H}\ominus BH^2_\cal{H}$ and these may fail to decompose as a direct sum of invariant subspaces contained in $H^2$. Therefore, one expects the scalar theory and vector-valued theory to be fundamentally different. A first indication of this fact is given by \cite{DPRS}*{Theorem 5.3} and Theorem~\ref{cstarthm} of this paper provides an extension of this result.
\section{A Distance Formula}\label{distance}
Let $\cal{A}$ be a \wk{}-closed, unital subalgebra of $H^\infty$, let $z_1,\ldots,z_n\in \bb{D}$ and let $\cal{I}$ denote the ideal 
\[\cal{I} := \{f\in \cal{A}\,:\, f(z_1)=\cdots = f(z_n) = 0\}.\] 
The ideal $\cal{I}$ is \wk{}-closed and the codimension of $\cal{I}$ in $\cal{A}$ is at most $n$.

In this section we give a formula for $\norm{f+\cal{I}}$, where $f\in L^\infty$. The formula relates the norm of $\norm{f+\cal{I}}$ to the norm of certain off-diagonal compressions of the operator $M_f$. The result is valid for subalgebras $A\subseteq H^\infty$ that have a property which we call predual factorization.

We identify $L^\infty$ as the dual of $L^1$ and refer to $L^1$ as the predual of $L^\infty$. 

We will say that a subspace $\cal{X}\subseteq H^\infty$ has \textit{predual factorization} if the following two properties hold
\begin{enumerate}
\item There exists a subspace (not necessarily closed) $\cal{S}\subseteq L^1$ with $[S]_1=\cal{X}_\perp$
\item Given $f\in \cal{S}$, there exists an inner function $\phi$ such that $\phi f\in H^1$.
\end{enumerate}

A simple consequence of Riesz factorization is that any function $f\in \cal{S}$ can be written as $f=\psi u^2$ where $\psi$ is unimodular, $u$ is outer and $\mod{f}^{1/2} = \mod{u}$. 

\begin{prop}
Suppose that $\{\cal{X}_j\,:\,j\in J\}$ is a set of \wk-closed subspaces of $H^\infty$. If for each $j\in J$, $\cal{X}_j$ has predual factorization, then $\cal{X}:=\bigcap_{j\in J}\cal{X}_j$ has predual factorization.
\end{prop}

\begin{proof}
We note that 
\[
\cal{X}_\perp = \left(\bigcap_{j\in J} \cal{X}_j\right)_\perp = \left[\bigcup_{j\in J}\{(\cal{X}_j)_\perp\,:\, j\in J\}\right]_1. 
\]
Set $\cal{S} = \linspan\{(\cal{X}_j)_\perp\,:\, j\in J\}$. Given $f_j\in \cal{X}_j$ there exists an inner function $\phi_j$ such that $\phi_j f_j\in H^1$. If $f=\sum_{i=1}^m c_i f_{j_i} \in \cal{S}$, then $\phi f \in H^1$ where $\phi = \phi_{j_1}\cdots \phi_{j_m}$. Hence, $\cal{S}$ has predual factorization.
\end{proof}

\begin{prop}
If $\cal{X}$ is a subspace of $L^\infty$ such that $BH^\infty \subseteq \cal{X}$, then $\cal{X}$ has predual factorization.
\end{prop}

\begin{proof}
We have $\cal{X}_\perp \subseteq (BH^\infty)_\perp = \cl{B}H^1_0$. The inner function $B$ multiplies $\cal{X}_\perp$ into $H^1$.
\end{proof}

\begin{cor}
If $\{B_j\,:\,j\in J\}$ is a set of inner functions, and $\cal{X}_j\supseteq B_jH^\infty$,  then $\cal{X}=\bigcap_{j\in J} \cal{X}_j$ has predual factorization.
\end{cor}

\begin{cor}
If $\{B_j\,:\,j\in J\}$ is a set of inner functions, then $\cal{A} = \bigcap_{j\in J} H^\infty_{B_j}$ has predual factorization.
\end{cor}

Recall that a function $u\in H^2$ is called \textit{outer} if $[H^\infty u] = H^2$. Given an outer function $u\in H^2$ we define $\cal{M}_u=[\cal{A}u]$, $\cal{K}_u$ to be the span of the kernel functions for $\cal{M}_u$ at the points $z_1,\ldots,z_n$ and $\cal{N}_u:=\cal{M}_u\ominus \cal{K}_u=\{f\in \cal{M}_u\,:\,f(z_1)=\cdots=f(z_n) = 0\}.$ Given a subspace $\cal{M}\subseteq L^2$ we denote by $P_\cal{M}$ the orthogonal projection of $L^2$ onto $\cal{M}$. 

\begin{lemma}\label{Ihasprefact}
let $z_1,\ldots,z_n$ be $n$ points in $\bb{D}$ and suppose $\cal{A}$ has predual factorization. If $\cal{I}$ is the ideal of functions in $\cal{A}$ such that $f(z_j)=0$, for $j=1,\ldots,n$, then $\cal{I}$ has predual factorization.
\end{lemma}

\begin{proof}
Let $\cal{I}_\perp$ be the preannihilator of $\cal{I}$ in $L^1$. Since $\cal{A}$ has predual factorization there exists a subspace $\cal{S}\subseteq \cal{A}_\perp$ such that 
\begin{enumerate}
\item The closure of $\cal{S}$ in the $L^1$ norm is $\cal{A}_\perp$
\item For each $f\in \cal{S}$, there exists an inner function $\phi$ such that $\phi f\in H^1$.
\end{enumerate} 
Note that $\cal{I}_\perp=\cal{A}_\perp+\linspan\{\overline{k_{z_j}}\,:\,1\leq j \leq n\}$, where $k_z$ is the Szeg\"o kernel at the point $z$. If $E$ is the Blaschke product for the points $z_1,\ldots,z_n$, then $E\overline{k_{z_j}}\in H^\infty$ for $j=1,\ldots,n$. The space $\tilde{\cal{S}}=\cal{S}+\linspan\{\overline{k_{z_j}}\,:\,1\leq j \leq n\}$ is dense in $\cal{I}_\perp$. Given $h+v\in \tilde{\cal{S}}$, with $h\in \cal{A}_\perp$ and $v\in \linspan\{\overline{k_{z_j}}\,:\,1\leq j \leq n\}$, there exists an inner function $\phi$ such that $\phi h\in H^1$ and so $E\phi(h+v)\in H^1$ with $E\phi$ is inner. Hence, $\cal{I}_\perp$ has predual factorization.
\end{proof}

\begin{lemma}
Let $\cal{I}$, $\cal{A}$ be as in Lemma~\ref{Ihasprefact}. If $u$ is an outer function, then $[\cal{I}u] = \cal{N}_u$.
\end{lemma}

\begin{proof}
Since every function in $\cal{I}$ vanishes at $z_1,\ldots,z_n$, $[\cal{I}u]\subseteq \cal{N}_u$. On the other hand given $f\in \cal{N}_u$ we know that there exists $f_m\in \cal{A}$ such that $\norm{f_m u -f}_2\to 0$. Since $u$ does not vanish at any point of the disk we see that $f_m(z_j)\to 0$ for $j=1,\ldots,n$. By a construction similar to the one in Lemma~\ref{pointsep} we see that there exists functions $e_j\in \cal{A}$ such that $e_j(z_i)=\delta_{i,j}$. Setting $g_m=f_m-\sum_{i=1}^n f_m(z_i)e_i$ we see that $g_m u$ converges to $f$ in $H^2$ and $g_m\in \cal{I}$. Hence, $\cal{N}_u\subseteq [\cal{I}u]$ and our proof is complete.
\end{proof}

We will now prove our distance formula. 

\index{Distance formula}
\begin{thm}\label{distanceLinftytoI}
Let $z_1,\ldots,z_n$ be $n$ points in $\bb{D}$, let $\cal{A}$ be a \wk-closed, unital subalgebra of $H^\infty$ with predual factorization and let $\cal{I}$ be the ideal of functions in $\cal{A}$ such that $f(z_j)=0$ for $j=1,\ldots,n$.  If $f\in L^\infty$, then 
\[
\norm{f+\cal{I}}=\sup_{u}\norm{(I-P_{\cal{N}_u})M_fP_{\cal{M}_u}},
\]
where the supremum is taken over all outer functions $u\in H^2$.
\end{thm}

\begin{proof}
We have,
\begin{align*}
\norm{f+\cal{I}}&=\sup\left\{\mod{\int fg}\,:\, g\in \cal{I}_\perp,\norm{g}_1\leq 1\right\}\\
&=\sup\left\{\mod{\int fg}\,:\, g\in \tilde{\cal{S}},\norm{g}_1\leq 1\right\},
\end{align*}
where $\tilde{S}$ is a dense subspace of $\cal{I}_\perp$ with the property that each function in $\cal{S}$ can be multiplied into $H^1$ by an inner function.
Let $g\in \tilde{\cal{S}}$, and let  $\phi$ be inner with the property that $\phi g \in H^1$ and factor $\phi g$ as $g_1 u$ where $g_1, u\in H^2$, $u$ is outer, and $\norm{u}_2=\norm{g_1}_2=\norm{g}_1^{1/2}$. It follows that $g=g_2u$ where $g_2\in L^2$ and $u$ is outer with $\norm{u}_2=\norm{g_2}_2=\norm{g}_1^{1/2}$. 

Since $g\in \cal{I}_\perp$, for all $h\in \cal{I}$ we get
\[
0 = \int gh = \int g_2uh =\inp{hu}{\cl{g_2}}. 
\]
This shows $\cl{g_2}\perp [\cal{I}u]=\cal{N}_u$. Hence 
\begin{align*}
\mod{\int fg}&=\mod{\inp{fu}{\cl{g_2}}}\\
& =\mod{\inp{f P_{\cal{M}_u}u}{(I-P_{\cal{N}_u})\cl{g_2}}}\\
&\leq \norm{(I-P_{\cal{N}_u})M_fP_{\cal{M}_u}}.
\end{align*}

For the other inequality we let $h\in \cal{I}$. We have $M_h\cal{M}_u \subseteq \cal{N}_u$ and so $(I-P_{\cal{N}_u})M_h P_{\cal{M}_u} = 0$. Therefore,
\begin{align*}
\norm{(I-P_{\cal{N}_u})M_fP_{\cal{M}_u}} & = \norm{(I-P_{\cal{N}_u})M_{f+h}P_{\cal{M}_u}} \\
&\leq \norm{M_{f+h}} \leq \norm{f+h}_\infty \leq \norm{f+\cal{I}}
\end{align*}
\end{proof}

We point out that the proof of Theorem~\ref{distanceLinftytoI} holds in the case $n=0$ to give the distance of an element in $L^\infty$ from the algebra $\cal{A}=\bigcap_{j\in J} H^\infty_{B_j}$. This result can be interpreted as a Nehari-type theorem for the algebra $\cal{A}$.

\begin{thm}\label{NehariforA}
If $f\in L^\infty$, then $\norm{f+\cal{A}}=\sup_{u}\norm{(I-P_{\cal{M}_u})M_fP_{\cal{M}_u}}$.
\end{thm}

\section{Interpolation}\label{interp}

Let $\cal{A}$ be a subalgebra of $H^\infty$, unital and $\wk$-closed as before. Given $n$ points $z_1,\ldots,z_n$ in the disk $\bb{D}$ and $n$ complex numbers $w_1,\ldots,w_n$, the interpolation problem for $\cal{A}$ is to determine conditions for the existence $f\in \cal{A}$ with $\norm{f}\leq 1$ such that $f(z_j)=w_j$ for $j=1,\ldots, n$. Such an $f$ will be called a solution. For the algebra $H^\infty$, the Nevanlinna-Pick theorem gives us a necessary and sufficient condition for the existence of a solution. 

Suppose that  $f_1,f_2\in \cal{A}$ such that $f_i(z_j)=w_j$ for $j=1,\ldots,n$, $i=1,2$. If $\cal{I}$ denotes the ideal of functions in $\cal{A}$ such that $f(z_j)=0$ for $j=1,\ldots,n$, then $f_1-f_2\in \cal{I}$. If we assume the existence of at least one function $f\in\cal{A}$ such that $f(z_j)=w_j$, then all other solutions are of the form $f+g$ for some $g\in \cal{I}$. As $\cal{A}$ is $\wk$-closed it follows that a solution exists if and only if $\norm{f+\cal{I}}\leq 1$. 

In this section we will prove the interpolation theorem for algebras with predual factorization. While the interpolation theorem is a little abstract we will see that it contains as special cases the original Nevanlinna-Pick theorem and the interpolation result from~\cite{DPRS}. It also provides us with an interpolation theorem for algebras of analytic functions on embedded disks.

\begin{thm}\label{Hinftybinterp1}
Let $\cal{A}$ be a \wk-closed, unital subalgebra of $H^\infty$ which has predual factorization. Let $z_1,\ldots,z_n\in \bb{D}$ and $w_1,\ldots,w_n\in \bb{C}$. Let $K^u$ denote the kernel function for the space $\cal{M}_u:=[\cal{A}u]$. There exists a function $f\in \cal{A}$ such that $f(z_j)=w_j$, $j=1,\ldots,n$ with $\norm{f}_\infty\leq 1$ if and only if for all outer functions $u\in H^2$, $[(1-w_i\overline{w_j})K^u(z_i,z_j)]_{i,j=1}^n\geq 0$.
\end{thm}

The proof of this result will follow from the distance formula obtained in Theorem~\ref{distanceLinftytoI}. 

We first need to establish the fact that the multiplier algebra of $\cal{M}_u$ contains $\cal{A}$. 

\index{Multiplier algebra of $H^\infty_B$}
\begin{prop}\label{multalg}
Let $\cal{A}$ be a \wk-closed, unital subalgebra of $H^\infty$ and let $u\in H^2$ be an outer function. If $\cal{M}:=[\cal{A}u]$, then $\cal{A}\subseteq \mult(\cal{M})$.
\end{prop}

\begin{proof}
It is straightforward that $\cal{A}(\cal{M})\subseteq \cal{M}$. Since $u$ does not vanish on the disk we see that none of the kernel functions in $\cal{M}$ are the zero function. If $M_f$ denotes the multiplication operator on $\cal{M}$ induced by $f$, then $\norm{f}_{\mult}=\norm{M_f}_{B(\cal{M})}\geq \norm{f}_\infty$. On the other hand if $h\in \cal{M}\subseteq L^2$, then
\[
\norm{M_fh}^2 = \int \mod{fh}^2\leq  \norm{f}_\infty^2\norm{h}^2, 
\]
which proves $\norm{f}_{\mult}\leq \norm{f}_\infty$.  
\end{proof}

In the special case where $\cal{H}=\bigcap_{j\in J} H^2_{B_j}$ is $\cal{A}=\bigcap_{j\in J} H^\infty_{B_j}$ we can improve on the previous proposition.

\begin{prop}\label{multalg2}
Let $\{B_j\,:\,j\in J\}$ be a set of inner functions. The multiplier algebra $\mult(\bigcap_{j\in J} H^2_{B_j}) = \bigcap_{j\in J} H^\infty_{B_j}$.
\end{prop}

\begin{proof}
Let us denote $\cal{M}:=\bigcap_{j\in J} H^2_{B_j}$. Let $f\in \mult(\cal{M})$. Since $1\in \cal{M}$ none of the kernel functions in $\cal{M}$ can be zero. This shows that any $f\in \mult(\cal{M})$ must be bounded. If $f\in \mult(\cal{M})$, then  $f\in \cal{M}$, since $1\in \cal{M}$. Hence, $f=\lambda_j + B_j k_j$ for $j\in J$, $k_j\in H^2$. Since $f$ is bounded so is $k_j$ and we have shown that $f\in \bigcap_{j\in J} H^\infty_{B_j}$. On the other hand any function $f\in \bigcap_{j\in J} H^\infty_{B_j}$ multiplies $\cal{M}$ into itself. It remains to be seen that $\norm{f}_{\mult} \leq \norm{f}_\infty$. This follows from 
\[
\norm{M_fh}^2 = \int \mod{fh}^2\leq  \norm{f}_\infty^2\norm{h}^2, 
\]
where $h\in \cal{M}\subseteq L^2$. 
\end{proof}

We now have all, but one, of the pieces required for our interpolation theorem. In order to use the distance formula to deduce the interpolation theorem we must know that there exists an interpolating function $f\in\cal{A}$. The proof of this fact is contained in Lemma~\ref{pointsep}. Let us assume the lemma for a moment and see how the interpolation theorem follows.  

If $f\in \cal{A}$, then $M_f$ leaves $\cal{M}_u$ invariant and 
\begin{align}\label{dist2}
\norm{f+\cal{I}}&=\sup_{u}\norm{(I-P_{\cal{N}_u})M_fP_{\cal{M}_u}}\\
&=\sup_{u}\norm{(I-P_{\cal{M}_u}+P_{\cal{K}_u})M_fP_{\cal{M}_u}}\\
&=\sup_{u}\norm{P_{\cal{K}_u}M_fP_{\cal{M}_u}}\\
\label{dist2f}&=\sup_{u}\norm{M_f^\ast P_{\cal{K}_u}}.
\end{align}

If $k^u_z$ denotes the kernel function for $\cal{M}_u$ at $z$, then a spanning set for $\cal{K}_u$ is given by $\{k^u_{z_1},\ldots,k^u_{z_n}\}$. Standard results about multiplier algebras of reproducing kernel Hilbert spaces tell us that the norm of $M_f^\ast P_{\cal{K}_u}$ is at most 1 if and only if $[(1-w_i\overline{w_j})K^u(z_i,z_j)]_{i,j=1}^n\geq 0$. Combining this fact with equation~\eqref{dist2}--\eqref{dist2f} proves the interpolation theorem.

Note, when $\cal{A} = H^\infty$, the cyclic subspace $[H^\infty u] = H^2$ and so we recover the classical Nevanlinna-Pick theorem. When $\cal{A}=H^\infty_B$ the structure of the cyclic subspaces is given by the next Proposition.

\begin{prop}\label{cyclic1}
If $B$ is an inner function, then $[H^\infty_B u] = [v]\oplus BH^2$, where $v = P_{H^2\ominus BH^2}u$. 
\end{prop}

\begin{proof}
Let $u=v\oplus Bw$. If $v=0$, then $u \in BH^2$ which contradicts the fact that $u$ is outer. If $\lambda+Bf \in H^\infty_B$, then $(\lambda+Bf)u = \lambda v+B(v+fv + Bfw) \in v\oplus BH^2$. Conversely, let $f=\lambda v\oplus Bh \in ([v]\oplus BH^2) \ominus [H^\infty_B u]$. Since $f\perp [H^\infty_B u]$ we see that $0=\inp{\lambda v+Bh}{Bgu}=\inp{h}{gu}$ for all $g\in H^\infty$. Hence, $h\perp [H^\infty u]=H^2$ and so $h=0$. Now, 
\[0=\inp{\lambda v}{u}=\lambda \inp{v}{v+Bw}=\lambda \inp{v}{v}=\lambda \norm{v}^2\]
and so $\lambda=0$. 
\end{proof}

Hence, the collection of cyclic subspaces for the algebra $H^\infty_B$ is a contained in the collection of subspaces of the form $[v]\oplus BH^2$. We can assume of course, that $v$ is a unit vector. The vector $v$ has an additional property. If $\phi$ is an inner function that divides $v$ and $B$, then $\phi$ divides $u$. Since $u$ is outer, we see $\gcd(v,B)=1$.

Let $v\in H^2\ominus BH^2$ be a unit vector and set $H^2_v=[v]\oplus BH^2$. This is a reproducing kernel Hilbert space with kernel function,
\[K^v(z,w) = v(z)\cl{v(w)}+\frac{B(z)\cl{B(w)}}{1-z\cl{w}}.\]
Let $\cal{K}_v=\linspan\{k_{z_1}^v,\ldots,k_{z_n}^v\}$ and let $\cal{N}_v=\linspan\{f\in H^2_v\,:\, f(z_j) = 0, j=1,\ldots,n\}$. Applying our distance formula to $H^\infty_B$ we get the following analogue of Nehari's theorem.
\begin{thm}
If $f\in L^\infty$, then 
\[\norm{f+\cal{I}}=\sup\norm{(1-P_{\cal{N}_v})M_fP_{H^2_v}},\]
where the supremum is taken over all unit vectors $v\in H^2\ominus BH^2$.
\end{thm}

The interpolation result now reads as follows:
\begin{thm}
Let $z_1,\ldots,z_n\in \bb{D}$ and $w_1,\ldots,w_n\in \bb{C}$. There exists a function $f\in H^\infty_B$ such that $f(z_j)=w_j$ if and only if the matrices,
\[\left[(1-w_i\cl{w_j})K^v(z_i,z_j)\right],\]
are positive for all unit vectors $v\in H^2\ominus BH^2$.
\end{thm}
Note in both the above theorems that the collection of subspaces is really parametrized by the one-dimensional subspaces of $H^2\ominus BH^2$.

The purpose of the next two lemmas is to show that if the matrix 
\[
\left[(1-w_i\overline{w_j})K^u(z_i,z_j)\right]_{i,j=1}^n 
\]
is positive for just one outer function $u\in H^2$, then there exists an interpolating function for the algebra $\cal{A}$. 

\begin{lemma}
Let $\cal{H}$ be a finite-dimensional Hilbert space, let $v_1,\ldots,v_n$ be a basis for $H$ and let $W_1,\ldots,W_n,W_{n+1}\in M_p$. Suppose that $v_{n+1}\in \cal{H}$ and $v_{n+1}=\sum_{i=1}^n \alpha_i v_i$.  If the matrix $Q=\left[(I-W_iW_j^\ast)\inp{v_j}{v_i}\right]_{i,j=1}^{n+1}$ is positive, then for all $1\leq i\leq n$ either $\alpha_i=0$ or $W_i=W_{n+1}$.
\end{lemma}

\begin{proof}
Let $W_s=\left[w^{(s)}_{i,j}\right]_{i,j=1}^p$ and consider the matrix $Q_k$ that we get by compressing to the $(k,k)$ entry of each block in $Q$. The $(i,j)$-th entry of this matrix is $\left(1-\sum_{l=1}^{p}w_{k,l}^{(i)}\overline{w_{k,l}^{(j)}}\right)\inp{v_j}{v_i}$. Let $\lambda_1,\ldots,\lambda_{n+1}\in \bb{C}$ and note that
\[\sum_{i,j=1}^{n+1}\left(1-\sum_{l=1}^{p}w_{k,l}^{(i)}\overline{w_{k,l}^{(j)}}\right)\inp{v_j}{v_i}\overline{\lambda_i}\lambda_j\geq 0.\]
By setting $\lambda_j=\alpha_j$ for $j=1,\ldots,n$ and $\lambda_{n+1}=-1$ we get,
\begin{align*}
0\leq &\sum_{i,j=1}^n\left(1-\sum_{l=1}^{p}w_{k,l}^{(i)}\overline{w_{k,l}^{(j)}}\right)\inp{v_j}{v_i}\overline{\alpha_i}\alpha_j\\
&-\sum_{i=1}^n \left(1-\sum_{l=1}^{p}w_{k,l}^{(i)}\overline{w_{k,l}^{(n+1)}}\right)\inp{v_{n+1}}{v_i}\overline{\alpha_i}\\
&-\sum_{j=1}^n \left(1-\sum_{l=1}^{p}w_{k,l}^{(n+1)}\overline{w_{k,l}^{(j)}}\right)\inp{v_j}{v_{n+1}}\alpha_j\\
&+\left(1-\sum_{l=1}^{p}\mod{w_{k,l}^{(n+1)}}^2\right)\norm{v_{n+1}}^2.
\end{align*}
This simplifies to 
\[\norm{v_{n+1}-\sum_{i=1}^n \alpha_iv_i}^2-\sum_{l=1}^{p}\norm{\sum_{i=1}^n(w_{k,l}^{(n+1)}-w_{k,l}^{(i)})\alpha_i v_i }^2\geq 0\]
which gives $\sum_{i=1}^n(w_{k,l}^{(n+1)}-w_{k,l}^{(i)})\alpha_i v_i=0$, for $1\leq k, l\leq p$. If $\alpha_i\not =0$, then by the linear independence of $v_1,\ldots,v_n$ we get $w_{k,l}^{(n+1)}=w_{k,l}^{(i)}$ and so $W_{n+1}=W_i$. 
\end{proof}

\begin{lemma}\label{pointsep}
Let $\cal{A}$ be a unital, \wk-closed subalgebra of $H^\infty$. Let $u$ be an outer function and let $\cal{M}=[\cal{A}u]$. Let $K$ be the kernel function of $\cal{M}$, $z_1,\ldots,z_n$ be $n$ points in the disk and $W_1,\ldots,W_n\in M_k$. If 
\begin{equation}\label{matpos}
\left[(1-W_iW_j^\ast)K(z_i,z_j)\right]_{i,j=1}^n\geq 0, 
\end{equation}
then there exists $F\in M_k(\cal{A})$ such that $F(z_j)=W_j$.
\end{lemma}

\begin{proof}
We may assume after reordering the points that $\{k_{z_1},\ldots.k_{z_m}\}$ is  basis of $\linspan\{k_{z_j}\,:\,1\leq j\leq n\}$, with $m\leq n$. Since $u$ is outer, $u$ is non-zero at every point of $\bb{D}$. There exists $f\in \cal{M}$ such that $\dfrac{f(z_j)}{u(z_i)}$ are distinct for for $1\leq j\leq m$. If this is not the case, then $\cl{u(z_j)}^{-1}k_{z_j}-\cl{u(z_i)}^{-1}k_{z_i} = 0$ is a non-trivial linear combination of $k_{z_i}$ and $k_{z_j}$. Since $\cal{M}$ is the closure of $[\cal{A}u]$, we conclude that there exists $g\in \cal{A}$ such that $g(z_j)$ are distinct for $1\leq j\leq m$.  
By setting 
\begin{equation}
\label{pointsep1}e_j=\prod_{r=1,r\not= j }^m\frac{g-g(z_r)}{g(z_j)-g(z_r)}\in \cal{A},
\end{equation}
we see that $e_i(z_j)=\delta_{i,j}$, for $1\leq i,j\leq m$. 
Let $h=\sum_{i=1}^m w_i e_i$ and note that for $1\leq j \leq m$, $h(z_j)=w_j$. To complete the proof we need to show that $h(z_j)=w_j$ for $j>m$.

Let $j>m$ and suppose that $k_{z_j} = \sum_{l=1}^m \alpha_l k_{z_l}$. We have seen that the matrix positivity condition~\eqref{matpos} implies that either  $w_j=w_l$ or $\alpha_l=0$, for $1\leq l\leq m$. Hence,
\begin{align*}
h(z_j)&=\sum_{i=1}^m w_ie_i(z_j) =\sum_{i=1}^m w_iu(z_j)^{-1}(ue_i)(z_j) \\
&=\sum_{i=1}^m w_i u(z_j)^{-1}\left(\sum_{l=1}^m \cl{\alpha_l} (ue_i)(z_l)\right)\\
&= \sum_{i=1}^m w_i u(z_j)^{-1}\cl{\alpha_i}u(z_i)= u(z_j)^{-1}\sum_{i=1}^m w_i \cl{\alpha_i} u(z_i) \\
&= w_ju(z_j)^{-1}\sum_{i=1}^m \cl{\alpha_i} u(z_i) = w_ju(z_j)^{-1}u(z_j)=w_j.
\end{align*}
The matrix case follows easily.
\end{proof}

\section{The $C^\ast$-envelope of $H^\infty_B/\cal{I}$.}\label{cstar}

The $C^\ast$-envelope of an operator algebra was defined by Arveson~\cite{arveson:cstar2}. Loosely speaking the $C^\ast$-envelope of an operator algebra $\cal{A}$, which is denoted $C^\ast_e(\cal{A})$, is the smallest $C^\ast$-algebra on which $\cal{A}$ has a completely isometric representation. Arveson's work established the existence of the $C^\ast$-envelope in the presence of what are called boundary representations. The existence of the $C^\ast$-envelope of an operator algebra was established in full generality by Hamana, whose approach did not have any relation to boundary representations.
 
\index{C@$C^\ast$-envelope}
\begin{thm}[Arveson-Hamana]\label{cstarenv}
Let $\cal{A}$ be an operator algebra. There exists a $C^\ast$-algebra, which is denoted $C^\ast_e(\cal{A})$, such that 
\begin{enumerate}
\item\label{cstarenv1} There is a completely isometric representation $\gamma:\cal{A} \to C^\ast_e(\cal{A})$.
\item\label{cstarenv2} Given a completely isometric representation $\sigma:\cal{A} \to \cal{B}$, where $\cal{B}$ is a $C^\ast$-algebra and $C^\ast(\sigma(\cal{A}))= \cal{B}$, there exists an onto $\ast$-homomorphism $\pi:\cal{B}\to C^\ast_e(\cal{A})$ such that $\pi \circ \sigma = \gamma$. 
\end{enumerate}
\end{thm}
It is easy to see that the $C^\ast$-envelope is essentially unique up to $\ast$-isomorphism. For a detailed description of the $C^\ast$-envelope we refer the reader to \cite{cbmoa}.

For the algebra $H^\infty$ the $C^\ast$-envelope of $H^\infty/\cal{I}$ is $M_n$. In \cite{So} the algebras 
\[
H^\infty_{a_1,\ldots,a_m}:=\{f\in H^\infty\,:\,f(a_1)=\ldots=f(a_m)\},
\]
were examined and the following result was obtained.

\begin{thm}[Sollazo \cite{So}]\label{solazzo}
Let $a_1=0$, $a_2=\frac{1}{2}$ and let $z_1=0$, $z_2,z_3\in \bb{D}$ with $z_1,z_2,z_3$ distinct. The $C^\ast$-envelope of the algebra $H^\infty_B/\cal{I}$ is $M_4$. 
\end{thm}

Note that the quotient $H^\infty_B/\cal{I}$ is a $3$-idempotent algebra \cite{idem}. When we compare this to the classical case we see there has been a jump in the dimension of the $C^\ast$-envelope from $3$ to $4$. 

This \index{Dimension jump}\textit{dimension jump} phenomenon has also been observed in \cite{DPRS}*{Theorem~5.3} for the algebra $\bb{C}+z^2H^\infty$. In this section we will show, given certain constraints on the number of zeros in the Blaschke product $B$, that a similar result is true for the algebra $H^\infty_B/\cal{I}$. The first step in understanding the quotient $H^\infty_B/\cal{I}$ is to gain some knowledge about the structure of the ideal $\cal{I}$. 

We will consider only the case where $B$ is a finite Blaschke product. To fix notation we let $\alpha_1,\ldots,\alpha_p$ be the zeros of $B$, we assume that these are distinct and have multiplicity $m_j\geq 1$ and we set $m=m_1+\ldots+m_p$. We arrange the points $z_1,\ldots,z_n$ so that $B(z_j)=0$ for $j=1,\ldots,r$ and $B(z_j)\not =0$ for $j=r+1,\ldots, n$. Denote by $E$ the Blaschke product for the points $z_1,\ldots,z_n$. It is clear that $\cal{I}=H^\infty_B\cap EH^\infty$. Since $B$ is a finite Blaschke product we see that $W\subseteq (H^2\ominus BH^2)\cap H^\infty$ and $\cal{I}=E(W+BH^\infty)$. This can also be seen directly from the fact that $\cal{I}$ is invariant for $H^\infty_B$.

\begin{thm}\label{idealstruct}
Let $B$ be a finite Blaschke product and let $\cal{I}$ be the ideal of functions in $H^\infty_B$ that vanish at the $n$ points $z_1,\ldots,z_n$. If $r=0$, then 
\[
\cal{I} = E([w]+BH^\infty)
\]
for some $w\in H^\infty\cap(H^2\ominus BH^2)$. If $r\geq 1$, then 
\[
\cal{I} = \lcm(B,E)H^\infty=E(W+BH^\infty),
\]
where $W$ is $r$-dimensional.
\end{thm}

\begin{proof}
Let $f\in \cal{I}$ and write $f=\lambda+Bg\in EH^\infty$, where $g\in H^\infty$. By evaluating at $z_1,\ldots,z_n$ we obtain $\lambda + B(z_j)g(z_j)=0$. 

First, consider the case where $r=0$. We can write $f=\lambda +B(\sum_{j=1}^n c_j k_{z_j})+BEh$ for some choice of $c_1,\ldots,c_n\in \bb{C}$ and $h\in H^\infty$. Hence, $\lambda +B(\sum_{j=1}^n c_j k_{z_j})$ is $0$ at the points $z_1,\ldots,z_n$ and so $\lambda +B(z_i) \sum_{j=1}^n c_jK(z_i,z_j)=0$ for $i=1,\ldots,n$. Rewriting this as a linear system we get 
\[
\begin{bmatrix} B(z_1) \\ & \ddots \\ & & B(z_n) \end{bmatrix}\begin{bmatrix}K(z_1,z_1) & \cdots & K(z_1,z_n) \\ \vdots & & \vdots \\ K(z_n,z_1) & \cdots & K(z_n,z_n)\end{bmatrix}\begin{bmatrix}c_1 \\ \vdots \\ c_n\end{bmatrix}=-\lambda \begin{bmatrix}1 \\ \vdots \\ 1\end{bmatrix}. 
\]
Since $r=0$, this system has a unique solution and the constants $c_1,\ldots,c_n$ can be taken to depend  linearly on $\lambda$. In this case $W$ is one-dimensional.

If $r\geq 1$, then $\lambda =0$ and $g(z_j)=0$ for $j=r+1,\ldots,n$. Hence, $f=B\phi_{z_{r+1}}\cdots \phi_{z_n}h$, $f\in \mathrm{lcm}(B,E)H^\infty$ and $\cal{I}\subseteq \mathrm{lcm}(B,E)H^\infty$. The reverse inclusion is straightforward. Let $C=B\cl{\gcd(B,E)}\in H^2$. From
\begin{align*}
\lcm(B,E)H^2 &= \lcm(B,E)\left([k_{z_1},\ldots,k_{z_r}]\oplus\gcd(B,E)H^2\right)\\
&=EB\cl{\gcd(B,E)}\left([k_{z_1},\ldots,k_{z_r}]\oplus\gcd(B,E)H^2\right)\\
&=E(C[k_{z_1},\ldots,k_{z_r}]\oplus BH^2),
\end{align*}
we see that $W$ is $r$-dimensional. 
\end{proof}

For an outer function $u$, 
\[
[\cal{I}u]=[\mathrm{lcm}(B,E)H^\infty u]=\mathrm{lcm}(B,E)[H^\infty u]=\mathrm{lcm}(B,E)H^2.  
\]
We have seen in Proposition~\ref{cyclic1} that $[H^\infty_B u]=[v]\oplus BH^2$ for some vector $v$ and so 
\begin{align*}
\cal{K}_u&=([v]\oplus BH^2)\ominus \mathrm{lcm}(B,E)H^2\\
&=[v]\oplus B(H^2\ominus \phi_{z_{r+1}}\cdots \phi_{z_{n}}H^2)\\
&=[v]\oplus B[k_{z_{r+1}},\ldots,k_{z_n}].
\end{align*}
The space $\cal{K}_u$ has dimension $(n-r)+1$. Note that this is also the dimension of the quotient algebra $H^\infty_B/\cal{I}$. Our distance formula says that interpolation is possible if and only if the compression of $M_f^\ast$ to $[v]\oplus B[k_{z_{r+1}},\ldots,k_{z_n}]$ is a contraction for all $v\in H^2\ominus BH^2$. 

In the case where one or more of the points $z_1,\ldots,z_n$ is a zero of $B$, i.e., when $r\geq 1$, the distance of $f\in L^\infty$ from $\cal{I}$ is the distance of $f$ from  $\mathrm{lcm}(B,E)H^\infty$. This is the case we will examine more closely. The objective will be to show that the scalar-valued result in Theorem~\ref{Hinftybinterp1} is not the correct matrix-valued interpolation result. This result also generalizes the result from~\cites{So}.

A basis for $\cal{K}=H^2\ominus \mathrm{lcm}(B,E)H^2$ is given by the vectors 
\begin{equation}\label{Kbasis}
\cal{E}:=\{z^i k_{\alpha_j}^{i+1}\,:\,1\leq j \leq p, 0\leq i \leq m_j-1\}\cup\{k_{z_{r+1}},\ldots,k_{z_n}\}.
\end{equation}
We begin by computing the matrix of $M_f^\ast|_\cal{K}$ with respect to the basis $\cal{E}$. It is an elementary calculation to show, for $f\in H^2$ and $m\geq 0$, that
\[
\frac{f^{(m)}(w)}{m!} = \inp{f}{z^mk^{m+1}_w}.
\]

\begin{lemma}\label{derivativebasis}
If $f\in H^\infty$, then 
\[M_f^\ast(z^m k_w^{m+1})=\sum_{j=0}^m\frac{1}{j!}\overline{f^{(j)}(w)}z^{m-j}k_w^{m-j+1}.\]
\end{lemma}

\begin{proof}
Let $g\in H^2$ and consider
\begin{align*}
\inp{g}{M_f^\ast z^mk_w^{m+1}}&=\inp{fg}{z^mk_w^{m+1}}=\frac{(fg)^{(m)}(w)}{m!}\\
&=\frac{1}{m!}\sum_{j=0}^m \binom{m}{j}f^{(j)}(w)g^{(m-j)}(w)\\
&=\frac{1}{m!}\sum_{j=0}^m \binom{m}{j}f^{(j)}(w)(m-j)!\inp{g}{z^{m-j}k_w^{m-j+1}}\\
&=\inp{g}{\sum_{j=0}^m\frac{1}{j!}\overline{f^{(j)}(w)}z^{m-j}k_w^{m-j+1}}.
\end{align*}
From this, we see that 
\[
M_f^\ast(z^mk_w^{m+1})=\sum_{j=0}^m\frac{1}{j!}\overline{f^{(j)}(w)}z^{m-j}k_w^{m-j+1}.
\]
\end{proof}

When $f\in H^\infty_B$, Lemma~\ref{derivativebasis} and the constraints on page~\pageref{constraints} show us that the matrix of $M_f^\ast$ is diagonal with respect to the basis $\cal{E}$. The matrix of $M_f^\ast|_\cal{K}$ is given by

\[
D_f^\ast = \begin{bmatrix} \overline{f(\alpha_1)}I_{m_1} \\  & \ddots \\ &  &  \overline{f(\alpha_p)}I_{m_p} \\ & & & \overline{f(z_{r+1})} \\ & & & & \ddots \\ & & & & & \overline{f(z_n)}  \end{bmatrix},
\]
If we partition the basis $\cal{E}$ as in~\eqref{Kbasis}, then the grammian matrix with respect to this basis has the form
\[
Q=\begin{bmatrix} Q_1 & Q_2 \\ Q_2^\ast & P \end{bmatrix},
\]
where $P$ is the Pick matrix for the points $z_{r+1},\ldots,z_n$. Since $Q$ is the grammian matrix of a linearly independent set it is invertible and positive. The matrix $Q_1$ is $m\times m$, positive and invertible, and the matrix $Q_2$ is an $m\times (n-r)$ matrix of rank $\min\{m,n-r\}$. 

For a function $f\in H^\infty$, Sarason's generalized interpolation~\cites{Sa} shows that the distance of $f$ from the ideal $\phi H^\infty$, i.e., $\norm{f+\phi H^\infty}$, is given by the norm of the compression of $M_f$ to $H^2\ominus \phi H^2$. This distance formula is also valid in the matrix-valued case. 

Let $T$ be an operator on a finite-dimensional Hilbert space $\cal{H}$, of dimension $n$ say, let $\cal{E}$ be a Hamel basis for $\cal{H}$ and let $A$ be the matrix of $T$ with respect to $\cal{E}$. The operator $T$ is a contraction if and only if $I_n - A^\ast QA\geq 0$, where $Q$ is the grammian matrix with respect to $\cal{E}$. 

Using this last fact, if $f\in H^\infty_B$, then 
\begin{align*}
\norm{M_f^\ast|_\cal{K}}\leq 1 & \iff  Q-D_f QD_f^\ast \geq 0\\ 
& \iff  Q^{1/2}(I-Q^{-1/2}D_f QD_f^\ast Q^{-1/2})Q^{1/2}\geq 0\\
& \iff I-Q^{-1/2}D_f QD_f^\ast Q^{-1/2}\geq 0 \\
& \iff I-(Q^{-1/2}D_f Q^{1/2})(Q^{-1/2}D_f Q^{1/2})^\ast\geq 0. 
\end{align*}
This induces a completely isometric embedding $\rho$ of $H^\infty_B/\cal{I}$ in $M_{m+n-r}$ given by
\[
\rho(f)=Q^{-1/2}D_f Q^{1/2}.
\]

The universal property of the $C^\ast$-envelope tells us that $C^\ast_e(H^\infty_B/\cal{I})$ is a quotient of $\cal{B}:=C^\ast(\rho(H^\infty_B/\cal{I}))$. Since we are dealing with a representation on a finite-dimensional space we know that $\cal{B}$ is a direct sum of matrix algebras. In the event that $\cal{B}=M_{m+n-r}$ we see that $\cal{B}=C^\ast_e(H^\infty_B/\cal{I})$. This follows from the fact that $M_{m+n-r}$ is simple. 

\begin{thm}
Let $r\geq 1$ and let $\cal{B}$ be the $C^\ast$-subalgebra of $M_{m+n-r}$ generated by the image of $\rho$. The algebra $\cal{B}=M_{m+n-r}$ if and only if $m\leq n-r$.
\end{thm}

\begin{proof}
We examine the commutant of $\cal{B}$ and show that $\cal{B}'$ contains only scalar multiples of the identity. Let $R=Q^{1/2}$ and let $RXR^{-1}\in \cal{B}'$. 

It is possible to choose $f\in H^\infty_B$ such that  $f(\alpha_i)=1$ for all $1\leq i \leq p$ and $f(z_j)=0$ for $r+1 \leq j \leq n$. Given $j$, with $r+1\leq j\leq n$, it is possible to choose $f$ such that $f(z_j)=1$, $f(\alpha_i)=f(z_l)=0$ for $1\leq i\leq p$ and $l\not= j$.  Therefore $\cal{B}$ is generated by $R^{-1}E_jR$ where $E_0:=E_{1,1}+\ldots+E_{m,m}$ and $E_j:=E_{m+j,m+j}$ for $1 \leq j \leq n-r$.

The matrix $RXR^{-1}\in \cal{B}'$ if and only if 
\[
RXR^{-1}R^{-1}E_j R = R^{-1}E_j R RXR^{-1}
\]
and 
\[
RXR^{-1}(R^{-1}E_j R)^\ast = (R^{-1}E_j R)^\ast RXR^{-1}.
\]
This happens if and only if $Q X Q^{-1} E_j= E_j QXQ^{-1}$ and $XE_j=E_jX$. These conditions tell us that $X$ and $QXQ^{-1}$ are both block diagonal with 1 block of size $m\times m$ followed by $n-r$ blocks of size $1$. Let us write 
\[
X=\begin{bmatrix} A & 0 \\ 0 & D\end{bmatrix},
QXQ^{-1}=\begin{bmatrix} B & 0 \\ 0 & E\end{bmatrix},
\]
where $D$ and $E$ are scalar diagonal of size $(n-r)$. We have,
\[
\begin{bmatrix} Q_1 & Q_2 \\ Q_2^\ast & P \end{bmatrix} \begin{bmatrix} A & 0 \\ 0 & D\end{bmatrix} = \begin{bmatrix} B & 0 \\ 0 & E\end{bmatrix} \begin{bmatrix} Q_1 & Q_2 \\ Q_2^\ast & P \end{bmatrix} .
\]
This tells us that $PD=EP$, where $P=[p_{i,j}]_{i,j=1}^{(n-r)}$ is the Pick matrix. Since $p_{i,j}d_i = p_{i,j}e_j$ and $p_{i,j}$ are non-zero for $1\leq i,j \leq (n-r)$, we get $d_i=e_j$ for $1\leq i,j\leq (n-r)$. Hence, we may assume that $D=E=I_{n-r}$. Now comparing the off-diagonal entry we see that $Q_2^\ast A = Q_2^\ast$, $BQ_2 = Q_2$ and so 
\[
Q_2^\ast A = Q_2^\ast B^\ast = Q_2^\ast.
\]
Rewriting this we get 
\begin{equation}\label{q2eq}Q_2^\ast(A-B^\ast)=Q^\ast_2(I_m -A)=Q^\ast_2(I_m -B^\ast).\end{equation}

If $m \leq n-r$, then $Q_2^\ast$ has rank $m$ which implies $I_m=A=B$, $X=I_{m+n-r}$ and 
\[
\cal{B}=\cal{B}''=\{I_{m+n-r}\}'=M_{m+n-r}.
\]
On the other hand if $m>n-r$, then there exist $m-n+r$ linearly independent solutions to the equation $Q_2^\ast v =0$. These can be used to construct matrices $A,B\not = I_m$ that solve equation (\ref{q2eq}). Hence, $\cal{B}\not = M_{m+n-r}$.
\end{proof}

\begin{thm}\label{cstarthm}
Let $r\geq 1$. The $C^\ast$-envelope of $H^\infty_B/\cal{I}$ is $M_{m+n-r}$ if and only if $m\leq n-r$. 
\end{thm}

\begin{proof}
This follows from Hamana's theorem and the fact that $M_{m+n-r}$ is simple.
\end{proof}

As a corollary we obtain the following generalization of a theorem from \cite{DPRS}.

\begin{cor}[\cite{DPRS}*{Theorem~5.3}]
Let $z_1=0$ and $n\geq 3$. The $C^\ast$-envelope of $H^\infty_1/\cal{I}$ is $M_{n+1}$.
\end{cor}

\begin{proof}
Since $r=1$ and $n\geq 3$ we see that $n-r = m = 2$. Hence, by the previous result $C^\ast_e(H^\infty_1/\cal{I}) = M_{m+n-r} = M_{n+1}$.
\end{proof}

As a corollary we also obtain Solazzo's result \cite{So}, which we stated as Theorem~\ref{solazzo}, for the algebra $H^\infty_{0,\frac{1}{2}}$. 

To close our discussion we want to make a few statements about the relevance of the result on $C^\ast$-envelopes in distinguishing between scalar-valued and matrix-valued problems.

Note that the collection of one-dimensional subspaces of $H^2\ominus BH^2$ can be identified with the complex-projective $m$-sphere $X=PS^m$. For a point $v\in X$ let us denote by $H^2_v$ the  subspace $[v]\oplus BH^2$ and let the kernel for $H^2_v$ be denoted $K^v$. For a fixed pair of points $z,w\in \bb{D}$, the map $(z,w)\mapsto K^v(z,w)$ is continuous. Denote by $\cal{K}_v$ the span of the kernel functions at the points $z_1,\ldots,z_n$ for $H^2_v$. The interpolation theorem tells us that there is an isometric representation of $H^\infty_B/\cal{I}$ on $C(X,M_{(n-r)+1})$ given by $\sigma(f+\cal{I})(v)=P_{\cal{K}_v}M_fP_{\cal{K}_v}$. If $\sigma$ is a completely isometric representation, then $\cal{C}=C^\ast(\sigma(H^\infty_B/\cal{I}))$ is a candidate for $C^\ast_e(H^\infty_B/\cal{I})$. However, the $C^\ast$-algebra $\cal{C}$ is a subalgebra of $M_{(n-r)+1}(C(X))$ and as such its irreducible representations can be at most $(n-r+1)$-dimensional. The fact that $m\geq 2$ tells us that $m+(n-r)>(n-r)+1$ and this implies that $C^\ast_e(H^\infty_B/\cal{I})$ cannot be contained completely isometrically in $M_{(n-r)+1}(C(X))$. This contradiction proves that the matrix-valued analogue of the interpolation result in Theorem~\ref{Hinftybinterp1} is generally false.

\end{document}